\documentclass[a4paper,11pt]{article}

\usepackage[margin=1.0in]{geometry}

\usepackage{bm}
\usepackage{graphicx}
\usepackage{tikz}
\usetikzlibrary{arrows.meta,decorations.pathreplacing}
\usepackage[bb=boondox]{mathalfa} 
\usepackage[mathlines]{lineno}

\usepackage{color}

\usepackage{amssymb}    
\usepackage{amsthm}     
\usepackage{thmtools}   
\usepackage{mathtools}  
\usepackage{mathrsfs}   
\usepackage{dsfont}     
\usepackage{csquotes}

\newcommand{\dist}{\operatorname{dist}}

\theoremstyle{plain}
\newtheorem{theorem}{Theorem}[section]
\newtheorem{corollary}[theorem]{Corollary}
\newtheorem{lemma}[theorem]{Lemma}
\newtheorem{proposition}[theorem]{Proposition}
\newtheorem{conjecture}[theorem]{Conjecture}

\theoremstyle{definition}

\newtheorem{example}{Example}
\newtheorem{observation}[theorem]{Observation}

\newcommand{\mc}{\mathcal}

\DeclareMathOperator*{\bigboxplus}{\tikz[baseline,line width=.2ex]{\node[minimum size=.9em,draw,anchor=base,inner sep=0ex]{\Large $+$};}\hspace{.05cm}}

\author
{Lorenzo Ciardo\thanks{Department of Computer Science, University of
Oxford, Wolfson Building, Parks Road, OX1 3QD Oxford, UK. Email: \texttt{lorenzo.ciardo@cs.ox.ac.uk}.}
}
\title{Perron value and moment of rooted trees}

\begin{document}

\maketitle

\begin{abstract}
\noindent The Perron value $\rho(T)$ of a rooted tree $T$ has a central role in the study of the algebraic connectivity and characteristic set, and it can be considered a weight of \textit{spectral} nature for $T$. A different, \textit{combinatorial} weight notion for $T$ -- the moment $\mu(T)$ -- emerges from the analysis of Kemeny's constant in the context of random walks on graphs. In the present work, we compare these two weight concepts showing that $\mu(T)$ is \enquote{almost} an upper bound for $\rho(T)$ and the ratio $\mu(T)/\rho(T)$ is unbounded but at most linear in the order of $T$. To achieve these primary goals, we introduce two new objects associated with $T$ -- the Perron entropy and the neckbottle matrix -- and we investigate how different operations on the set of rooted trees affect the Perron value and the moment.
\end{abstract}





\section{Introduction} 
\label{sec_introduction}
A major factor for the success of spectral graph theory as an approach to study graph properties is its ability to meaningfully quantify how connected a graph is and provide center notions. The Laplacian matrix $L$ of a graph $G=(V(G),E(G))$ -- defined by $L=D-A$, where $D$ and $A$ are the diagonal degree matrix and adjacency matrix of $G$, respectively -- is particularly well suited to capture the concepts of connectivity and centrality. Its second-smallest eigenvalue, which we shall denote by $a(G)$, is positive if and only if $G$ is connected and does not decrease when a new edge is inserted into $G$. For these reasons, it was named \textit{algebraic connectivity} by Fiedler \cite{Fiedler_alg_conn} and, since then, its properties have been extensively investigated in relation to combinatorial connectivity notions such as the vertex connectivity, the edge connectivity, and the number of cut vertices \cite{Fiedler_alg_conn,Kirkland_a_bound_cutpoints,Kirkland_an_upper_bound_many_cutpoints}, the isoperimetric number \cite{Alon,Mohar_isoperimetric}, the genus \cite{Molitierno_genus}, and other graph-theoretic parameters \cite{Fallat_Kirkland}. Additionally, the \textit{Fiedler vectors} of $G$ -- i.e., the eigenvectors for $L$ corresponding to $a(G)$ -- can be used to identify a set of central vertices in $G$ known as the \textit{characteristic vertices} \cite{Fiedler_eigenvectors}.

Henceforth, we shall suppose that $G$ is a tree. In this case, the algebraic connectivity and the characteristic vertices are closely linked to the so-called \textit{bottleneck matrices} associated with rooted subtrees of $G$. The principal submatrix $L_v$ of $L$ obtained by removing from $L$ the row and column corresponding to a vertex $v$ is invertible, and its inverse $L_v^{-1}$ is a block diagonal matrix where each block $M=[m_{ij}]$ corresponds to a connected component (branch) $T$ of $G-v$. The entries of $M$ have a simple combinatorial description: $m_{ij}$ is the number of vertices in $T$ that simultaneously lie in the path joining $i$ to $r$ and in the path joining $j$ to $r$, where $r$ is the vertex of $T$ adjacent to $v$ in the original tree $G$ \cite{Kirkland_Neumann_Shader}. In particular, it is clear from this description that the entries of $M$ only depend on $T$ and on the vertex $r$ adjacent to $v$ in $G$, and are independent of the structure of $G-T$. As a consequence, $M$ is unambiguously associated with the pair $(T,r)$; we shall refer to it as to the bottleneck matrix of the \textit{rooted tree} $(T,r)$. Since $M$ is entrywise positive, by virtue of the Perron-Frobenius theory it has a simple positive dominant eigenvalue known as the \textit{Perron value} and denoted by $\rho((T,r))$ or simply $\rho(T)$ when the root $r$ is clear from the context (we also denote by $\rho(X)$ the spectral radius of a generic square matrix $X$). A corresponding eigenvector is called a \textit{Perron vector}. Finally, a \textit{Perron branch} for $G$ at $v$ is a branch attaining the maximum Perron value among all the branches of $G$ at $v$. The next result provides the link between the algebraic connectivity and characteristic vertices of a tree and the Perron value of its rooted subtrees\footnote{Analogous results can be obtained for generic connected graphs, by considering cut vertices and blocks instead of vertices and edges \cite{Kirkland_Fallat_Perron_components}.}. 
\begin{theorem}[\cite{Kirkland_Neumann_Shader}]
\label{FiedlerTheorem_Trees}
Let $G$ be an unrooted tree with more than one vertex. Exactly one of two cases occurs.
\begin{enumerate}
\item
There exists exactly one vertex $z$ such that $G$ has $k\geq 2$ Perron branches $B_1,B_2,\ldots,B_k$ at $z$.
$G$ is said to be a \textit{type I} tree and $z$ is its characteristic vertex. Moreover, in this case,
\begin{align*}
a(G)=\frac{1}{\rho(B_i)}\hspace{1.5cm} (i=1,2,\dots,k).
\end{align*}
\item
There exists exactly one edge $pq$ such that the unique Perron branch $B_p$ at $p$ contains $q$ and the unique Perron branch $B_q$ at $q$ contains $p$.
$G$ is said to be a \textit{type II} tree and $p,q$ are its characteristic vertices. Moreover, in this case,
\begin{align*}
a(G)=\frac{1}{\rho\left(M_p-\beta J\right)}=\frac{1}{\rho\left(M_q-(1-\beta)J\right)},
\end{align*}
where $M_p$ (resp. $M_q$) is the bottleneck matrix of $B_p$ (resp. $B_q$), $J$ is the all-ones matrix of suitable dimensions, and $\beta$ is a real number such that $0 < \beta < 1$.
\end{enumerate}
\end{theorem}
The fact that connectivity and center notions for trees can be determined by means of the Perron value of rooted trees provides the motivation to look for explicit expressions for this parameter or lower and upper bounds on its value. The former approach yielded exact formulae in the cases of rooted stars and paths \cite{Abreu_et_alii,EAGD17}:
\begin{align}
\label{expr_rho_star_2008}
\rho(\mathcal{S}_n)&=\frac{1}{2}\left(n+1+\sqrt{n^2+2n-3}\right)\!,\\
\label{expr_rho_path_2008}
\rho(\mathcal{P}_n)&=\frac{1}{2}\left(1-\cos\left(\frac{\pi}{2n+1}\right)\right)^{-1}\!\!\!,
\end{align}
where $\mathcal S_n$ denotes the rooted star on $n$ vertices with the central vertex as the root and $\mathcal{P}_n$ denotes the rooted path on $n$ vertices with one of the endpoints as the root.
The latter approach led to the definition of the so-called \textit{combinatorial Perron parameters} \cite{EALCGD19,EAGD17} -- lower bounds on the Perron value defined as the Rayleigh quotients of the bottleneck matrix and combinatorial surrogates of the Perron vectors. These parameters are shown to be close approximations of the Perron value, and their calculation is significantly faster since they do not require the computation of eigenvalues \cite{EALCGD19}.  

The theory of Markov chains provides a connectivity notion for graphs that is alternative to the algebraic connectivity. As shown in \cite{KemenySnell}, multiplying the \textit{mean first passage matrix} and the \textit{stationary distribution vector} of an irreducible discrete-time Markov chain results in a constant vector. The common value of the entries of this vector -- decreased by $1$ -- is known as \textit{Kemeny's constant} and expresses the expected transition time between two different states of the system, both randomly sampled according to the stationary distribution \cite{Levene_Loizou}. For the case of the random walk on an undirected graph, the corresponding Kemeny's constant provides a measure of the long-run ability of the graph structure to transmit information along its edges; hence, it can be seen as a connectivity notion for the graph. In \cite{Kemeny_bounds}, the authors study this parameter in the context of the random walk on a tree and show that it can be expressed in terms of Kemeny's constant for the random walks on certain subtrees by means of a quantity called \textit{moment}. The moment $\mu((T,r))$ (or simply $\mu(T)$) of a rooted tree $(T,r)$ is defined by
\begin{align}
\label{def_moment}
\mu((T,r))=\sum_{v\in V(T)}\dist(v,r)\deg(v),
\end{align} 
where $\dist(v,r)$ is the number of edges in the path joining $v$ to $r$ and $\deg(v)$ is the degree of $v$. This quantity is reminiscent of the homonymous notion in mechanics, once we consider the degree of a vertex to be proportional to its mass. Interpreting it as a \textit{combinatorial} weight for rooted trees, it looks natural to investigate its connection with a different weight of \textit{spectral} nature -- the Perron value. The goal of the current work is to explore how these two notions are related. 

Section \ref{sec_Perron_value_and_moment} contains the main results: we give upper (Theorem \ref{mu_rho_f}) and lower (Theorem \ref{mu_rho_4_7}) bounds for the Perron value in terms of the moment, and we show that the ratio $\mu(T)/\rho(T)$ is unbounded (Theorem \ref{thm_mu_rho_unbounded}). In the other sections, we develop the machinery useful to prove the results in Section \ref{sec_Perron_value_and_moment}. In particular, in Section \ref{sec_entropy}, the notion of \textit{Perron entropy} of rooted trees is introduced, as a measure of uniformity for the entries of the Perron vectors. In Section \ref{sec_neckbottle}, the so-called \textit{neckbottle matrix} -- closely related to the bottleneck matrix -- is considered. In Section  \ref{sec_rooted_sums_product_powers}, the behaviours of the Perron value and the moment are studied in connection with three different operations for rooted trees; as a by-product, we obtain a new lower bound for the algebraic connectivity of Bethe trees (Observation \ref{lower_bound_alg_conn_bethe}). Besides their use in this analysis, the notions of Perron entropy and neckbottle matrix can be of interest for future work on the Perron value of rooted trees and the algebraic connectivity.

\medskip 

\noindent\underline{Notation}:
We let $\mathbb{R}^n$ denote the space of $n$-dimensional real column vectors, and we identify such vectors with the corresponding $n$-tuples. The $i$'th standard unit vector is denoted by $e_i$, and the all-ones vector is denoted by $e$. To keep the notation light, for both $e_i$ and $e$ we do not indicate the dimension explicitly: it will be clear from the context. The Euclidean norm of a vector $w$ is denoted by $\|w\|$. The $n\times n$ identity matrix and all-ones matrix are denoted by $I_n$ and $J_n$, respectively. The set of nonnegative integers is denoted by $\mathbb{N}$, while the set of positive integers is denoted by $\mathbb{N}_{>0}$. The \textit{order} of a graph is the number of its vertices. The \textit{trivial} rooted tree -- denoted by $\mathcal{E}$ -- is the rooted tree of order $1$; a rooted tree is \textit{nontrivial} if its order is at least $2$.

\section{Perron entropy of rooted trees}
\label{sec_entropy}
Let $T$ be a rooted tree of order $n$, let $M$ be its bottleneck matrix, and let $w$ be a Perron vector for $M$. The quantity
\begin{align*}
H(T)=\frac{(e^Tw)^2}{\|w\|^2}
\end{align*}
is well defined since, by virtue of the Perron-Frobenius theory, the dimension of the eigenspace for $M$ corresponding to the Perron value is $1$; we shall refer to it as to the \textit{Perron entropy} of $T$. In Proposition \ref{bounds_perron_value_rooted_product}, this parameter will be used to express a lower bound on the Perron value of a particular product of rooted trees.
By applying the Cauchy–Schwarz inequality and the $1$-norm $2$-norm inequality, we see that
\begin{align*}
1\leq H(T)\leq n.
\end{align*}
The minimum value is attained only if $w$ is a multiple of a standard unit vector, while the maximum value is attained only if $w$ is a constant vector. Hence, we can interpret the Perron entropy as a measure of uniformity for the entries of the Perron vectors; this justifies the name chosen for the parameter.

Computational experiments show that, in most cases, the Perron entropy of a rooted tree decreases as the diameter increases, if the number of vertices is kept fixed. This accords with the fact -- noted in \cite{EAGD17} -- that a positive Perron vector can be approximated by a vector containing the distance of each vertex from the root. However, there are exceptions to this general trend. For example,
for the rooted path $\mathcal{P}_9$ and the rooted tree $T'$ obtained from $\mathcal{P}_8$ by attaching one additional pendent vertex to the vertex adjacent to the root, we compute $H(\mathcal{P}_9)\approx 7.665 > 7.660\approx H(T')$. In the remaining part of this section, we give explicit expressions for the Perron entropy of the rooted path $\mc P_n$ (Proposition \ref{entropy_path}) and the rooted star $\mc S_n$ (Proposition \ref{prop_entropy_star_2008}).
\begin{proposition}
\label{perron_vector_path}
For $n\in\mathbb{N}_{>0}$, each Perron vector of the bottleneck matrix of the rooted path $\mathcal{P}_n$ is a multiple of the vector $w=(w_i)$ defined by
\begin{align*}
w_i=\sin\left(\frac{i\pi}{2n+1}\right)\hspace{1.5cm}(i=1,2,\dots,n).
\end{align*} 
\end{proposition}
\begin{proof}
Let $U$ be the path on $2n+1$ vertices, considered as an unrooted tree. According to \cite[\S~1.4.4]{Brouwer_Haemers}, a Fiedler vector of $U$ is given by $z=(z_i)\in\mathbb{R}^{2n+1}$, where
\begin{align*}
z_i=\cos\left(\frac{(2i-1)\pi}{4n+2}\right)\hspace{1.5cm}(i=1,2,\dots,2n+1).
\end{align*}
Since $U$ is symmetric about the central vertex $v$, it is a type $I$ tree according to the classification in Theorem \ref{FiedlerTheorem_Trees}. Using \cite[Theorem~6.2.15]{Molitierno_book}, we have that the vector $(z_n,z_{n-1},\dots,z_1)$ is a Perron vector for the bottleneck matrix of $\mathcal{P}_n$. Observe that the $i$'th component of this vector is
\begin{align*}
z_{n-i+1}&=
\cos\left(\frac{(2(n-i+1)-1)\pi}{4n+2}\right)=
\cos\left(\frac{\pi}{2}-\frac{i\pi}{2n+1}\right)
=
\sin\left(\frac{i\pi}{2n+1}\right)
\end{align*}
as wanted.
\end{proof}
\noindent We shall use the following result from \cite{Knapp}.
\begin{theorem}[\cite{Knapp}]
\label{trigonometric_formulas}
Let $d$ be a real number that is not an integer multiple of $2\pi$ and let $N\in\mathbb{N}$. Then
\begin{align*}
\sum_{i=0}^N\sin(id)&=\frac{\sin\left((N+1)d/2\right)}{\sin(d/2)}\sin\left(Nd/2\right)
\intertext{and}
\sum_{i=0}^N\cos(id)&=\frac{\sin\left((N+1)d/2\right)}{\sin(d/2)}\cos\left(Nd/2\right)\!.
\end{align*}
\end{theorem}

\begin{proposition}
\label{entropy_path}
For $n\in\mathbb{N}_{>0}$, the Perron entropy of the rooted path $\mathcal{P}_n$ is
\begin{align*}
H(\mathcal{P}_n)=\frac{\cot^2\left(\frac{\pi}{4n+2}\right)}{2n+1}.
\end{align*}
\end{proposition}
\begin{proof}
For a positive integer $N$, using Theorem \ref{trigonometric_formulas}, we obtain
\begin{align*}
\sum_{i=1}^N\sin\left(\frac{i\pi}{N}\right)&=
\sum_{i=0}^N\sin\left(\frac{i\pi}{N}\right)=
\frac{\sin\left(\frac{(N+1)\pi}{2N}\right)}{\sin\left(\frac{\pi}{2N}\right)}\sin\left(\frac{\pi}{2}\right)\\
&=\frac{\sin\left(\frac{\pi}{2}+\frac{\pi}{2N}\right)}{\sin\left(\frac{\pi}{2N}\right)}=\frac{\cos\left(\frac{\pi}{2N}\right)}{\sin\left(\frac{\pi}{2N}\right)}=\cot\left(\frac{\pi}{2N}\right)\!.
\end{align*}
Hence, for the Perron vector $w$ of the bottleneck matrix of $\mathcal{P}_n$ given in Proposition \ref{perron_vector_path}, we have 
\begin{align*}
2e^Tw&=2\sum_{i=1}^n\sin\left(\frac{i\pi}{2n+1}\right)
=
\sum_{i=1}^n\sin\left(\frac{i\pi}{2n+1}\right)+
\sum_{i=1}^n\sin\left(\pi-\frac{i\pi}{2n+1}\right)\\
&=\sum_{i=1}^n\sin\left(\frac{i\pi}{2n+1}\right)+
\sum_{i=1}^n\sin\left(\frac{(2n+1-i)\pi}{2n+1}\right)\\
&=\sum_{i=1}^n\sin\left(\frac{i\pi}{2n+1}\right)+
\sum_{i=n+1}^{2n}\sin\left(\frac{i\pi}{2n+1}\right)\\
&=\sum_{i=1}^{2n}\sin\left(\frac{i\pi}{2n+1}\right)
=\sum_{i=1}^{2n+1}\sin\left(\frac{i\pi}{2n+1}\right)=\cot\left(\frac{\pi}{4n+2}\right)\!.
\end{align*}
Also, observe that
\begin{align*}
\|w\|^2&=\sum_{i=1}^nw_i^2=\sum_{i=1}^n \sin^2\left(\frac{i\pi}{2n+1}\right)=\sum_{i=1}^n\frac{1-\cos\left(\frac{2i\pi}{2n+1}\right)}{2}\\
&=\frac{n}{2}-\frac{1}{2}\sum_{i=1}^n\cos\left(\frac{2i\pi}{2n+1}\right)=
\frac{n+1}{2}-\frac{1}{2}\sum_{i=0}^n\cos\left(\frac{2i\pi}{2n+1}\right)\\
&=\frac{n+1}{2}-\frac{\sin\left(\frac{(n+1)\pi}{2n+1}\right)}{2\sin\left(\frac{\pi}{2n+1}\right)}\cos\left(\frac{n\pi}{2n+1}\right)=\frac{n+1}{2}-\frac{\sin\left(\frac{n\pi}{2n+1}\right)\cos\left(\frac{n\pi}{2n+1}\right)}{2\sin\left(\frac{\pi}{2n+1}\right)}\\
&=\frac{n+1}{2}-\frac{\sin\left(\frac{2n\pi}{2n+1}\right)}{4\sin\left(\frac{\pi}{2n+1}\right)}=\frac{n+1}{2}-\frac{\sin\left(\frac{\pi}{2n+1}\right)}{4\sin\left(\frac{\pi}{2n+1}\right)}=\frac{n+1}{2}-\frac{1}{4}=\frac{n}{2}+\frac{1}{4}.
\end{align*}
We conclude that
\begin{align*}
H(\mathcal{P}_n)&=\frac{(e^Tw)^2}{\|w\|^2}=\frac{\cot^2\left(\frac{\pi}{4n+2}\right)}{4\left(\frac{n}{2}+\frac{1}{4}\right)}
=
\frac{\cot^2\left(\frac{\pi}{4n+2}\right)}{2n+1}
\end{align*}
as desired.
\end{proof}
\begin{proposition}
\label{prop_entropy_star_2008}
For $n\in\mathbb{N}_{>0}$, the Perron entropy of the rooted star $\mc S_n$ is
\begin{align}
\label{expr_entropy_star_nontrivial}
H(\mc S_n)=\frac{(n^2+2n)\sqrt{n^2+2n-3}+n^3+3n^2-2}{(n+2)\sqrt{n^2+2n-3}+n^2+3n}
\end{align}
if $n\geq 2$, and $H(\mc S_n)=1$ if $n=1$.
\end{proposition}
\begin{proof}
The case $n=1$ is trivial, so we suppose $n\geq 2$. Let $r$ be the root of $\mc S_n$ and notice that its bottleneck matrix is
\begin{align*}
M=ee^T+I_n-e_re_r^T.
\end{align*}
From \eqref{expr_rho_star_2008}, we observe that the Perron value $\rho$ of $\mathcal{S}_n$ satisfies $\rho^2=\rho n+\rho-1$. The vector $w=\rho e-e_r$ is a Perron vector for $M$ since
\begin{align*}
Mw&=(ee^T+I_n-e_re_r^T)(\rho e-e_r)=\rho n e+\rho e-\rho e_r-e-e_r+e_r\\
&=(\rho n+\rho-1)e-\rho e_r=\rho^2 e-\rho e_r=\rho w.
\end{align*}
Therefore, the Perron entropy of $\mc S_n$ is
\begin{align}
\label{expr_H_star_partial}
H(\mc S_n)=\frac{(e^Tw)^2}{\|w\|^2}=\frac{\rho^2n^2-2\rho n+1}{\rho^2 n-2\rho +1}.
\end{align}
Plugging the formula \eqref{expr_rho_star_2008} into \eqref{expr_H_star_partial} yields the final expression \eqref{expr_entropy_star_nontrivial}. 
\end{proof}
\section{Neckbottle matrix}
\label{sec_neckbottle}
Let $T$ be a rooted tree having root $r$, whose vertices are labelled by $1,2,\dots,n$. We consider the partial order \enquote{$\preceq$} on the vertex set $V(T)$ defined by setting $j\preceq i$ if and only if the path joining $j$ to $r$ contains $i$ ($i,j\in V(T))$. From \cite[Lemma~2.1]{EAGD17}, we can express the bottleneck matrix $M$ of $T$ as $M=N^TN$, where $N$ is the \textit{path matrix} of $T$ -- i.e., the $n\times n$ $(0,1)$-matrix whose $(i,j)$'th entry is $1$ if $j\preceq i$, $0$ otherwise.
In this section we consider the matrix $Q=NN^T$, which we call the \textit{neckbottle matrix} of $T$ (the name suggests that the order of the two matrices $N^T$ and $N$ has been changed). From this definition, it follows that the $(i,j)$'th entry of $Q$ is the number of vertices $k$ such that $k\preceq i$ and $k\preceq j$.

Note that $\rho(M)=\rho(Q)$ since $M$ and $Q$ have the same eigenvalues. Therefore, the neckbottle matrix provides a new tool for finding or estimating the Perron value of a rooted tree.
\begin{observation}
\label{neckbottle_not_positive_1004_1236}
Unlike the bottleneck matrix, the neckbottle matrix can have zero entries. For example, the neckbottle matrix of the rooted star $\mathcal{S}_n$ is
\begin{align*}
Q=\begin{bmatrix}
n & e^T \\ 
e & I_{n-1}
\end{bmatrix}\!\!. 
\end{align*}
As a consequence, it can be convenient to compute the Perron value of a rooted tree by using the neckbottle matrix instead of the bottleneck matrix.
\end{observation}
The inverse of the bottleneck matrix $M$ has a combinatorial description: as mentioned in the Introduction, $M^{-1}$ is a specific submatrix of the Laplacian matrix of an unrooted tree containing $T$ as a subtree. We now show that the inverses of the path matrix and of the neckbottle matrix have a combinatorial description, too. For $i,j\in V(T)$, the expression $i\sim j$ indicates that $i$ is adjacent to $j$. Also, we say that $i$ and $j$ are \textit{brothers} if $i\neq j$, $\dist(i,r)=\dist(j,r)$, and there exists $k\in V(T)$ such that $i\sim k$ and $j\sim k$.
\begin{proposition}
\label{inverse_path_matrix}
Let $T$ be a rooted tree and let $N$ be its path matrix. Then $N$ is invertible and its inverse $X=[x_{ij}]$ satisfies
\begin{align}
\label{exp_inverse_path_matrix}
x_{ij}=\left\{
\begin{array}{cl}
1 & \mbox{if }\,  i=j\\
-1 & \mbox{if }\  i\sim j\,, j\preceq i\\
0 & \mbox{otherwise.} 
\end{array}
\right.
\end{align} 
\end{proposition} 
\begin{proof}
Multiplying the matrix $X$ defined in \eqref{exp_inverse_path_matrix} by $N$ yields
\begin{align*}
(XN)_{ij}=\sum_{k=1}^n x_{ik}n_{kj}=n_{ij}-\sum_{\substack{k\sim i\\k\preceq i}}n_{kj}=
\left\{
\begin{array}{cl}
1 & \mbox{if }\,i=j\\
0 & \mbox{if }\,j\not\preceq i\\
1-1=0  & \mbox{if }\, j\preceq i,\, i\neq j
\end{array}
\right.
\end{align*}
so $XN=I_n$ as desired.
\end{proof}
\begin{proposition}
Let $T$ be a rooted tree having root $r$ and let $Q$ be its neckbottle matrix. Then $Q$ is invertible and its inverse $Y=[y_{ij}]$ satisfies
\begin{align}
\label{expr_inverse_neckbottle_1354_0904}
y_{ij}=\left\{
\begin{array}{cl}
1 & \mbox{if }\, i=j=r\\
2 & \mbox{if }\, i=j\neq r\\
-1 & \mbox{if }\, i\sim j\\
1 & \mbox{if }\, i,j \mbox{ are brothers}\\
0 & \mbox{otherwise}.
\end{array}
\right.
\end{align} 
\end{proposition}
\begin{proof}
Since $Q=NN^T$, we have that $Q$ is invertible and
\begin{align*}
(Q^{-1})_{ij}&=((N^{-1})^TN^{-1})_{ij}=\sum_{k=1}^n(N^{-1})_{ki}(N^{-1})_{kj}.
\end{align*}
If $i=j$, we can write this as
\begin{align*}
(Q^{-1})_{ij}&=\sum_{k=1}^n(N^{-1})_{ki}^2=(N^{-1})_{ii}^2+\sum_{k\neq i}(N^{-1})_{ki}^2\\
&=1+|\{k\,:\, k\sim i,\, i\preceq k\}|
=
\left\{
\begin{array}{cl}
1 & \mbox{if }\, i=r\\
2 & \mbox{if }\, i\neq r.
\end{array}
\right.
\end{align*}
If $i\neq j$, we obtain
\begin{align*}
(Q^{-1})_{ij}&=(N^{-1})_{ii}(N^{-1})_{ij}+(N^{-1})_{ji}(N^{-1})_{jj}+\sum_{k\neq i,j}(N^{-1})_{ki}(N^{-1})_{kj}\\
&=
(N^{-1})_{ij}+(N^{-1})_{ji}+\sum_{k\neq i,j}(N^{-1})_{ki}(N^{-1})_{kj}.
\end{align*}
The last three lines of \eqref{expr_inverse_neckbottle_1354_0904} follow by observing that
\begin{align*}
(N^{-1})_{ij}+(N^{-1})_{ji}
=
\left\{
\begin{array}{cl}
-1 & \mbox{if }\, i\sim j\\
0 & \mbox{otherwise}
\end{array}
\right.
\end{align*}
and
\begin{align*}
\sum_{k\neq i,j}(N^{-1})_{ki}(N^{-1})_{kj}
&=
\left\{
\begin{array}{cl}
1 & \mbox{if }\,i,j \mbox{ are brothers}\\
0 & \mbox{otherwise.}
\end{array}
\right.\quad\qedhere
\end{align*}
\end{proof}
\begin{figure}
\centering
\begin{tikzpicture}[scale=1.25]
%
\draw [fill,draw=black] (-.1,-.1) rectangle ++(.2,.2);
\draw(0,.25) node{\small{$1$}};
\draw[fill]  (-1,-1) circle (0.1cm);
\draw(-1-.25,-1) node{\small{$2$}};
\draw[fill]  (1,-1) circle (0.1cm);
\draw(1+.25,-1) node{\small{$3$}};
\draw[fill]  (-1,-2) circle (0.1cm);
\draw(-1-.25,-2) node{\small{$4$}};
\draw[fill] (0.5,-2) circle (0.1cm);
\draw(0.5+.25,-2) node{\small{$5$}};
\draw[fill] (1.5,-2) circle (0.1cm);
\draw(1.5+.25,-2) node{\small{$6$}};
%
\draw[black] (0,0) -- (-1,-1) -- (-1,-2);
\draw[black] (0,0) -- (1,-1);
\draw[black] (1,-1) -- (0.5,-2);
\draw[black] (1,-1) -- (1.5,-2);
\end{tikzpicture}
\caption{A rooted tree. The root is indicated by a square.}
\label{fig_tree_examples_neckbottle}
\end{figure}
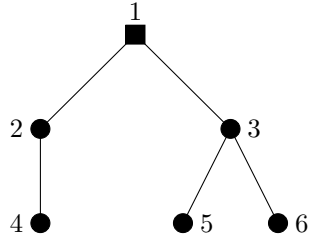
\begin{example}
For the rooted tree in Figure \ref{fig_tree_examples_neckbottle}, we find
\begin{align*}
N&=
\left[\begin{array}{cccccc} 1 & 1 & 1 & 1 & 1 & 1\\ 0 & 1 & 0 & 1 & 0 & 0\\ 0 & 0 & 1 & 0 & 1 & 1\\ 0 & 0 & 0 & 1 & 0 & 0\\ 0 & 0 & 0 & 0 & 1 & 0\\ 0 & 0 & 0 & 0 & 0 & 1 \end{array}\right]\!\!,
&
N^{-1}&=
\left[\begin{array}{cccccc} 1 & -1 & -1 & 0 & 0 & 0\\ 0 & 1 & 0 & -1 & 0 & 0\\ 0 & 0 & 1 & 0 & -1 & -1\\ 0 & 0 & 0 & 1 & 0 & 0\\ 0 & 0 & 0 & 0 & 1 & 0\\ 0 & 0 & 0 & 0 & 0 & 1 \end{array}\right]\!\!,
\\[5pt]
M&=
\left[\begin{array}{cccccc} 1 & 1 & 1 & 1 & 1 & 1\\ 1 & 2 & 1 & 2 & 1 & 1\\ 1 & 1 & 2 & 1 & 2 & 2\\ 1 & 2 & 1 & 3 & 1 & 1\\ 1 & 1 & 2 & 1 & 3 & 2\\ 1 & 1 & 2 & 1 & 2 & 3 \end{array}\right]\!\!,
&
M^{-1}&=
\left[\begin{array}{cccccc} 3 & -1 & -1 & 0 & 0 & 0\\ -1 & 2 & 0 & -1 & 0 & 0\\ -1 & 0 & 3 & 0 & -1 & -1\\ 0 & -1 & 0 & 1 & 0 & 0\\ 0 & 0 & -1 & 0 & 1 & 0\\ 0 & 0 & -1 & 0 & 0 & 1 \end{array}\right]\!\!,
\\[5pt]
Q&=
\left[\begin{array}{cccccc} 6 & 2 & 3 & 1 & 1 & 1\\ 2 & 2 & 0 & 1 & 0 & 0\\ 3 & 0 & 3 & 0 & 1 & 1\\ 1 & 1 & 0 & 1 & 0 & 0\\ 1 & 0 & 1 & 0 & 1 & 0\\ 1 & 0 & 1 & 0 & 0 & 1 \end{array}\right]\!\!,
&
Q^{-1}&=
\left[\begin{array}{cccccc} 1 & -1 & -1 & 0 & 0 & 0\\ -1 & 2 & 1 & -1 & 0 & 0\\ -1 & 1 & 2 & 0 & -1 & -1\\ 0 & -1 & 0 & 2 & 0 & 0\\ 0 & 0 & -1 & 0 & 2 & 1\\ 0 & 0 & -1 & 0 & 1 & 2 \end{array}\right]\!\!.
\end{align*}
\end{example}
\section{Rooted sum, product, and power}
\label{sec_rooted_sums_product_powers}
In this section, we consider three different operations on the set of rooted trees -- namely, the rooted sum, product, and power -- and we investigate how they affect the Perron value and the moment. The next, matrix-theoretic proposition, which is a direct consequence of a result in \cite[\S~3.5]{topics_matrix_analysis}, will be used in the analysis.
\begin{figure}
\centering
\begin{tikzpicture}[scale=.736]
\begin{scope}[shift={(0,0)}]
\draw  (.5,-0.2) node{\small{$T_1$}};
\draw [fill,draw=black] (-.1,-.1) rectangle ++(.2,.2);
\draw[fill] (-.5,-1) circle (0.1cm);
\draw[fill] (.5,-1) circle (0.1cm);
\draw[black] (-.5,-1) -- (0,0) -- (.5,-1);
\end{scope}
\begin{scope}[shift={(0,-3)}]
\draw  (.4,-.7) node{\small{$T_2$}};
\draw [fill,draw=black] (-.1,-.1) rectangle ++(.2,.2);
\draw[fill] (0,-1) circle (0.1cm);
\draw[fill] (0,-2) circle (0.1cm);
\draw[black] (0,0) -- (0,-1) -- (0,-2);
\end{scope}
\begin{scope}[shift={(0,-7)}]
\draw  (0.2,0.4) node{\small{$T_3$}};
\draw [fill,draw=black] (-.1,-.1) rectangle ++(.2,.2);
\end{scope}
\begin{scope}[shift={(3,0)}]
\draw  (1.5,-2.6) node{$a)\hspace{.3cm}\mbox{rooted sum}\hspace{.3cm} \displaystyle\bigboxplus_{i=1}^3T_i$};
\draw [fill,draw=black] (1.5+-.1,1+-.1) rectangle ++(.2,.2);
\draw[fill] (0,0) circle (0.1cm);
\draw[fill] (-.5,-1) circle (0.1cm);
\draw[fill] (.5,-1) circle (0.1cm);
\draw[fill] (1.5,0) circle (0.1cm);
\draw[fill] (1.5,-1) circle (0.1cm);
\draw[fill] (1.5,-2) circle (0.1cm);
\draw[fill] (3,0) circle (0.1cm);
\draw[black] (0,0) -- (1.5,1);
\draw[black] (1.5,0) -- (1.5,1);
\draw[black] (3,0) -- (1.5,1);
\draw[black] (-.5,-1) -- (0,0) -- (.5,-1);
\draw[black] (1.5,0) -- (1.5,-1) -- (1.5,-2);
\end{scope}
\begin{scope}[shift={(3,-5)}]
\draw  (1.5,-2.6) node{$b)\hspace{.3cm}\mbox{rooted product}\hspace{.3cm} \displaystyle T_1\boxtimes T_2$};
\draw[fill] (0,0) circle (0.1cm);
\draw[fill] (0,-1) circle (0.1cm);
\draw[fill] (0,-2) circle (0.1cm);
\draw [fill,draw=black] (1.5-.1,1-.1) rectangle ++(.2,.2);
\draw[fill] (1.5,1-1) circle (0.1cm);
\draw[fill] (1.5,1-2) circle (0.1cm);
\draw[fill] (3,0) circle (0.1cm);
\draw[fill] (3,-1) circle (0.1cm);
\draw[fill] (3,-2) circle (0.1cm);
\draw[black] (0,-2) -- (0,-1) -- (0,0) -- (1.5,1) -- (1.5,0) -- (1.5,-1);
\draw[black] (1.5,1) -- (3,0) -- (3,-1) -- (3,-2);
\end{scope}
\begin{scope}[shift={(12,-2)}]
\draw  (0,-5.1) node{$c)\hspace{.3cm}\mbox{rooted power}\hspace{.3cm} \displaystyle T_1^{\boxtimes 3}$};
\draw [fill,draw=black] (-.1,-.1) rectangle ++(.2,.2);
\draw[fill] (-.5,1) circle (0.1cm);
\draw[fill] (.5,1) circle (0.1cm);
\draw[fill] (-1.5,1) circle (0.1cm);
\draw[fill] (1.5,1) circle (0.1cm);
\draw[fill] (-1.5-0.5547,1+0.9707) circle (0.1cm);
\draw[fill] (-1.5-1.1094,1+0.1387) circle (0.1cm);
\draw[fill] (1.5+0.5547,1+0.9707) circle (0.1cm);
\draw[fill] (1.5+1.1094,1+0.1387) circle (0.1cm);
\draw[black](-.5,1)--(0,0)--(.5,1);
\draw[black](-1.5,1)--(0,0)--(1.5,1);
\draw[black](-1.5-0.5547,1+0.9707)--(-1.5,1)--(-1.5-1.1094,1+0.1387);
\draw[black](1.5+0.5547,1+0.9707)--(+1.5,1)--(1.5+1.1094,1+0.1387);
\begin{scope}[shift={(-.5*2,-0.866*2)},rotate=120]
\draw[fill] (0,0) circle (0.1cm);
\draw[fill] (-.5,1) circle (0.1cm);
\draw[fill] (.5,1) circle (0.1cm);
\draw[fill] (-1.5,1) circle (0.1cm);
\draw[fill] (1.5,1) circle (0.1cm);
\draw[fill] (-1.5-0.5547,1+0.9707) circle (0.1cm);
\draw[fill] (-1.5-1.1094,1+0.1387) circle (0.1cm);
\draw[fill] (1.5+0.5547,1+0.9707) circle (0.1cm);
\draw[fill] (1.5+1.1094,1+0.1387) circle (0.1cm);
\draw[black](-.5,1)--(0,0)--(.5,1);
\draw[black](-1.5,1)--(0,0)--(1.5,1);
\draw[black](-1.5-0.5547,1+0.9707)--(-1.5,1)--(-1.5-1.1094,1+0.1387);
\draw[black](1.5+0.5547,1+0.9707)--(+1.5,1)--(1.5+1.1094,1+0.1387);
\end{scope}
\begin{scope}[shift={(.5*2,-0.866*2)},rotate=2400]
\draw[fill] (0,0) circle (0.1cm);
\draw[fill] (-.5,1) circle (0.1cm);
\draw[fill] (.5,1) circle (0.1cm);
\draw[fill] (-1.5,1) circle (0.1cm);
\draw[fill] (1.5,1) circle (0.1cm);
\draw[fill] (-1.5-0.5547,1+0.9707) circle (0.1cm);
\draw[fill] (-1.5-1.1094,1+0.1387) circle (0.1cm);
\draw[fill] (1.5+0.5547,1+0.9707) circle (0.1cm);
\draw[fill] (1.5+1.1094,1+0.1387) circle (0.1cm);
\draw[black](-.5,1)--(0,0)--(.5,1);
\draw[black](-1.5,1)--(0,0)--(1.5,1);
\draw[black](-1.5-0.5547,1+0.9707)--(-1.5,1)--(-1.5-1.1094,1+0.1387);
\draw[black](1.5+0.5547,1+0.9707)--(+1.5,1)--(1.5+1.1094,1+0.1387);
\end{scope}
%
\draw[black](-.5*2,-0.866*2)--(0,0)--(+.5*2,-0.866*2);

\end{scope}
\end{tikzpicture}
\caption{Examples of rooted sum, product, and power for the three rooted trees on the left. The roots are indicated by squares.}
\label{fig_examples_operations_sum_prod_pow}
\end{figure}
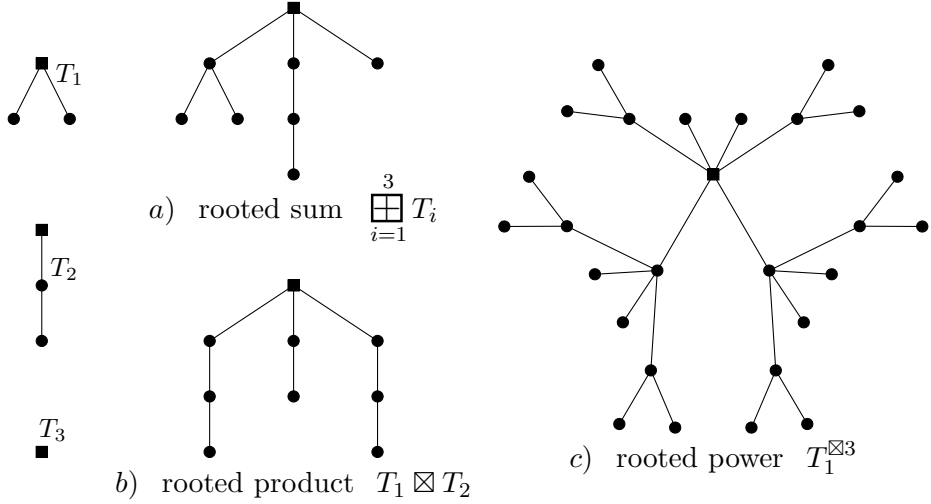
\begin{proposition}[\cite{topics_matrix_analysis}]
\label{prop_fundamental_topics}
Let $A=
\begin{bmatrix}
A_{11} & A_{12} \\ 
A_{12}^T & A_{22}
\end{bmatrix} 
$
be a symmetric positive semidefinite block matrix. Then $\rho(A)\leq \rho(A_{11})+\rho(A_{22})$.
\end{proposition}
Let $k\in \mathbb{N}_{>0}$. Given $k$ rooted trees $T_1,T_2,\dots,T_k$ having roots $r_i$ and orders $n_i$ $(i=1,2,\dots,k)$, we let their \textit{rooted sum} $\bigboxplus_{i=1}^k T_i$ be the rooted tree obtained by joining $r_1,r_2,\dots,r_k$ to an additional vertex $r$, which we take as the root (Figure \ref{fig_examples_operations_sum_prod_pow}~$a$). Observe that the order of $\bigboxplus_{i=1}^k T_i$ is $\sum_{i=1}^kn_i+1$.
\begin{proposition}
\label{neckbottle_composition}
Let $T_1,T_2,\dots,T_k$ be rooted trees and let $n=\sum_{i=1}^kn_i+1$. Then
\begin{align}
\label{ineq_rho_sum_18_08_2020}
\max_{1\leq i\leq k}\rho(T_i)\leq
\rho\left(\bigboxplus_{i=1}^k T_i\right)\leq \max_{1\leq i\leq k}\rho(T_i)+n.
\end{align}
\end{proposition}
\begin{proof}
Let $Q_i$ be the neckbottle matrix of $T_i$ $(i=1,2,\dots,k)$. Then, the neckbottle matrix of $\bigboxplus_{i=1}^k T_i$ is permutationally similar to
\begin{align*}
\begin{bmatrix}
n & x_1^T & x_2^T & \cdots & x_k^T \\ 
x_1 & Q_1 & O & \cdots & O \\ 
x_2 & O & Q_2 & \cdots & O \\ 
\vdots & \vdots & \vdots & \ddots & \vdots \\ 
x_k & O & O & \cdots & Q_k
\end{bmatrix} 
\end{align*}
for suitable vectors $x_1,x_2,\dots,x_k$, where $O$ denotes the zero block of suitable size. The first and second inequalities in \eqref{ineq_rho_sum_18_08_2020} follow from the Cauchy's interlacing theorem \cite[Theorem~4.3.28]{HJ} and Proposition \ref{prop_fundamental_topics}, respectively.
\end{proof}
\begin{proposition}
\label{moment_rooted_sum}
Let $T_1,T_2,\dots,T_k$ be rooted trees and let $n=\sum_{i=1}^kn_i+1$. Then
\begin{align*}
\mu\left(\bigboxplus_{i=1}^kT_i\right)=\sum_{i=1}^k\mu(T_i)+2n-2-k.
\end{align*}
\end{proposition}
\begin{proof}
In this proof, the expression $\deg_U$ (resp. $\dist_U$) shall indicate that the degree (resp. distance) is considered in the tree $U$. Denote $\bigboxplus_{i=1}^kT_i$ by $T$ and let $\vartheta=\sum_{i=1}^k\deg_{T_i}(r_i)$. We observe that
\begin{align*}
\mu(T)&=\sum_{v\in V(T)}\dist_T(v,r)\deg_T(v)=\sum_{i=1}^k\sum_{v\in V(T_i)}\dist_T(v,r)
\deg_T(v)\\
&=\sum_{i=1}^k \left(\dist_T(r_i,r)\deg_T(r_i)+
\sum_{\substack{v\in V(T_i)\\v\neq r_i}}\dist_T(v,r)
\deg_T(v)
\right)\\
&=\sum_{i=1}^k\left(
\deg_{T_i}(r_i)+1+
\sum_{\substack{v\in V(T_i)\\v\neq r_i}}
(\dist_{T_i}(v,r_i)+1)\deg_{T_i}(v)
\right)\\
&=
\vartheta+k+\sum_{i=1}^k\left(\mu(T_i)+2(n_i-1)-\deg_{T_i}(r_i)\right)\\
&=
\vartheta+k+\sum_{i=1}^k\mu(T_i)+2n-2-2k-\vartheta
=\sum_{i=1}^k\mu(T_i)+2n-2-k.\quad\qedhere
\end{align*}
\end{proof}
\begin{example}
\label{example_star_as_rooted_sum}
The rooted star may be expressed as the rooted sum of trivial trees: $\mathcal{S}_n=\bigboxplus_{i=1}^{n-1}\mc{E}$. Then, Proposition \ref{neckbottle_composition} yields the bound
$
\rho(\mathcal{S}_n)\leq 1+n
$,
which -- by virtue of the expression \eqref{expr_rho_star_2008} -- is asymptotically sharp as $n$ approaches infinity, while Proposition \ref{moment_rooted_sum} provides the exact value
$
\mu(\mathcal{S}_n)=n-1
$
as found in \cite{Kemeny_bounds}. 
\end{example}
Proposition \ref{neckbottle_composition} and Proposition \ref{moment_rooted_sum} allow to obtain results on the Perron value and the moment of a class of rooted trees that will prove useful in Section \ref{sec_Perron_value_and_moment}. For two integers $p\geq 1$ and $k\geq 2$, the \textit{rooted Bethe tree} $\mathcal B_{p,k}$ \cite{Heilmann_Lieb} is the rooted tree recursively defined as follows:
\begin{equation}
\label{def_bethe_trees}
\begin{array}{ll}
\mathcal{B}_{1,k}=\mathcal{E};&\\
\displaystyle\mc{B}_{p,k}=\bigboxplus_{i=1}^k \mc{B}_{p-1,k}&\hspace{.4cm}(p\geq 2).
\end{array}
\end{equation}
An example is shown in Figure \ref{fig_rooted_bethe_tree}. Observe that the order of $\mc{B}_{p,k}$ is 
\begin{align}
\label{order_bethe_trees}
|V(\mc{B}_{p,k})|=
\sum_{i=0}^{p-1}k^i=\frac{k^p-1}{k-1}.
\end{align}
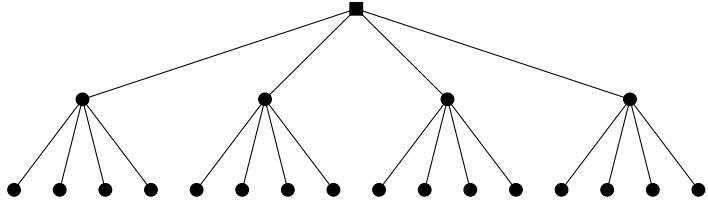
\begin{figure}
\centering
\begin{tikzpicture}[scale=1.2]
%
\draw [fill,draw=black] (-.07,-.07) rectangle ++(.14,.14);
\draw[fill]  (-3,-1) circle (0.07cm);
\draw[fill]  (-1,-1) circle (0.07cm);
\draw[fill]  (1,-1) circle (0.07cm);
\draw[fill]  (3,-1) circle (0.07cm);
\draw[fill]  (.25+.5*0,-2) circle (0.07cm);
\draw[fill]  (.25+.5*1,-2) circle (0.07cm);
\draw[fill]  (.25+.5*2,-2) circle (0.07cm);
\draw[fill]  (.25+.5*3,-2) circle (0.07cm);
\draw[fill]  (.25+.5*4,-2) circle (0.07cm);
\draw[fill]  (.25+.5*5,-2) circle (0.07cm);
\draw[fill]  (.25+.5*6,-2) circle (0.07cm);
\draw[fill]  (.25+.5*7,-2) circle (0.07cm);
\draw[fill]  (-.25-.5*0,-2) circle (0.07cm);
\draw[fill]  (-.25-.5*1,-2) circle (0.07cm);
\draw[fill]  (-.25-.5*2,-2) circle (0.07cm);
\draw[fill]  (-.25-.5*3,-2) circle (0.07cm);
\draw[fill]  (-.25-.5*4,-2) circle (0.07cm);
\draw[fill]  (-.25-.5*5,-2) circle (0.07cm);
\draw[fill]  (-.25-.5*6,-2) circle (0.07cm);
\draw[fill]  (-.25-.5*7,-2) circle (0.07cm);
%
\draw[black](0,0)--(-3,-1);
\draw[black](0,0)--(-1,-1);
\draw[black](0,0)--(1,-1);
\draw[black](0,0)--(3,-1);
\draw[black](-3,-1)--(-.25-.5*4,-2);
\draw[black](-3,-1)--(-.25-.5*5,-2);
\draw[black](-3,-1)--(-.25-.5*6,-2);
\draw[black](-3,-1)--(-.25-.5*7,-2);
\draw[black](-1,-1)--(-.25-.5*0,-2);
\draw[black](-1,-1)--(-.25-.5*1,-2);
\draw[black](-1,-1)--(-.25-.5*2,-2);
\draw[black](-1,-1)--(-.25-.5*3,-2);
\draw[black](1,-1)--(.25+.5*0,-2);
\draw[black](1,-1)--(.25+.5*1,-2);
\draw[black](1,-1)--(.25+.5*2,-2);
\draw[black](1,-1)--(.25+.5*3,-2);
\draw[black](3,-1)--(.25+.5*4,-2);
\draw[black](3,-1)--(.25+.5*5,-2);
\draw[black](3,-1)--(.25+.5*6,-2);
\draw[black](3,-1)--(.25+.5*7,-2);
\end{tikzpicture}
\caption{The rooted Bethe tree $\mathcal{B}_{3,4}$. The root is indicated by a square.}
\label{fig_rooted_bethe_tree}
\end{figure}
\begin{proposition}
\label{rho_ennary_tree}
Let $p\geq 1$ and $k\geq 2$ be integers. Then
\begin{align}
\label{eq_rho_ennary_tree}
\rho(\mc B_{p,k})\leq \frac{k^{p+1}-pk-k+p}{(k-1)^2}.
\end{align}
\end{proposition}
\begin{proof}
We use induction on $p$. If $p=1$, $\mc B_{p,k}=\mc E$ and both the left-hand side and the right-hand side of \eqref{eq_rho_ennary_tree} equal $1$. If $p\geq 2$, using Proposition \ref{neckbottle_composition} and the inductive hypothesis, we obtain
\begin{align*}
\rho(\mc B_{p,k})&=\rho\left(\bigboxplus_{i=1}^k\mc B_{p-1,k}\right)\leq \rho(\mc B_{p-1,k})+\frac{k^p-1}{k-1}\\
&\leq\frac{k^p-(p-1)k-k+(p-1)}{(k-1)^2}+ \frac{k^p-1}{k-1}=\frac{k^{p+1}-pk-k+p}{(k-1)^2}
\end{align*}
as desired.
\end{proof}
\begin{observation}
\label{lower_bound_alg_conn_bethe}
The result in Proposition \ref{rho_ennary_tree} can be used to provide a lower bound for the algebraic connectivity of an (unrooted) Bethe tree -- which, for the sake of simplicity, we shall indicate by the same notation as for its rooted counterpart. Observe that $\mc{B}_{p,k}$ is symmetric about the root. Hence, if $p\geq 2$, $\mc{B}_{p,k}$ is a type $I$ tree according to the classification in Theorem \ref{FiedlerTheorem_Trees}, where the unique characteristic vertex is the vertex $r$ corresponding to the root. Moreover, each branch at $r$ is a Perron branch isomorphic to $\mc{B}_{p-1,k}$. We conclude that
\begin{align}
\label{bound_alg_conn_bethe_tree}
a(\mc{B}_{p,k})=\frac{1}{\rho(\mc{B}_{p-1,k})}\geq 
\frac{(k-1)^2}{k^p-(p-1)k-k+(p-1)}
=\frac{(k-1)^2}{k^p-pk+p-1}.
\end{align}
Computational experiments show that the bound in \eqref{bound_alg_conn_bethe_tree} is quite tight. As an example, for the tree $\mc{B}_{6,6}$ (of order $9331$), the right-hand side of \eqref{bound_alg_conn_bethe_tree} produces a value that is the $99.96\%$ of the exact value of the algebraic connectivity. A different lower bound for the algebraic connectivity of Bethe trees was found in \cite{Rojo_Medina_tight} using matrix-theoretic techniques. We also point out that explicit formulae for some of the simple Laplacian eigenvalues of Bethe trees are known, see \cite{Rojo_Robbiano_explicit}. However, by virtue of \cite[Theorem~6.2.18]{Molitierno_book} (see also \cite[Theorem~2]{Grone_Merris}), the multiplicity of the algebraic connectivity as an eigenvalue of the Laplacian matrix of a type $I$ tree equals the number of Perron branches at the characteristic vertex decreased by $1$. Consequently, $a(\mc{B}_{p,k})$ is a simple eigenvalue of the Laplacian matrix of $\mc{B}_{p,k}$ if and only if $k=2$, so the formulae in \cite{Rojo_Robbiano_explicit} cannot be applied for $k\geq 3$.
\end{observation}
\begin{proposition}
\label{moment_ennary_tree}
Let $p\geq 1$ and $k\geq 2$ be integers. Then
\begin{align}
\label{expr_moment_bethe_tree}
\mu(\mc B_{p,k})=\frac{2pk^{p+1}-3k^{p+1}-2pk^{p}+k^{p}+k^2+k}{(k-1)^2}.
\end{align}
\end{proposition}
\begin{proof}
We use induction on $p$. If $p=1$, both the left-hand side and the right-hand side of \eqref{expr_moment_bethe_tree} equal $0$. If $p\geq 2$, using Proposition \ref{moment_rooted_sum} and the inductive hypothesis, we find 
\begin{align*}
\mu(\mc B_{p,k})&=\mu\left(\bigboxplus_{i=1}^k \mc{B}_{p-1,k}\right)=k\mu(\mc{B}_{p-1,k})+2\frac{k^p-1}{k-1}-2-k\\
&=k\frac{2(p-1)k^{p}-3k^{p}-2(p-1)k^{p-1}+k^{p-1}+k^2+k}{(k-1)^2}\\[-2pt]
&\quad+2\frac{k^p-1}{k-1}-2-k\\
&=\frac{2pk^{p+1}-3k^{p+1}-2pk^{p}+k^{p}+k^2+k}{(k-1)^2},
\end{align*}
thus validating the inductive step and concluding the proof.
\end{proof}

Let $T_1$ and $T_2$ be rooted trees having roots $r_i$ and orders $n_i$ $(i=1,2)$. Let $U_1,U_2,\dots,U_{n_1}$ be $n_1$ disjoint copies of $T_2$ and, for each edge $ij\in E(T_1)$, connect the root of $U_i$ to the root of $U_j$ with an edge. Finally, let the root of the tree thus constructed be the root of $U_{r_1}$. We denote the rooted tree obtained in this way by $T_1\boxtimes T_2$, and, following \cite{Godsil_McKay}, we refer to it as to the \textit{rooted product} of $T_1$ and $T_2$ (Figure \ref{fig_examples_operations_sum_prod_pow}~$b$). Observe that the order of $T_1\boxtimes T_2$ is $n_1n_2$ and, in general, $T_1\boxtimes T_2\neq T_2\boxtimes T_1$. Henceforth, we denote by $N_i$ (resp. $M_i$, $Q_i$) the path matrix (resp. bottleneck matrix, neckbottle matrix) of $T_i$ $(i=1,2)$. The \textit{Kronecker product} of an $m\times n$ matrix $A=[a_{ij}]$ and an $m'\times n'$ matrix $A'$ is the $mm'\times nn'$ block matrix
\begin{align*}
A\otimes A'=\begin{bmatrix}
a_{11}A' & a_{12}A' & \cdots & a_{1n}A' \\ 
a_{21}A' & a_{22}A' & \cdots & a_{2n}A' \\ 
\vdots & \vdots & \ddots & \vdots \\ 
a_{m1}A' & a_{m2}A' & \cdots & a_{mn}A'
\end{bmatrix}\!\!. 
\end{align*}
\begin{proposition}
\label{expr_path_matrix_rooted_product}
Let $T_1$ and $T_2$ be rooted trees. Then the path matrix of $T_1\boxtimes T_2$ is permutationally similar to
\begin{align*}
I_{n_1}\otimes N_2 +(N_1-I_{n_1})\otimes e_{r_2}e^T.
\end{align*}
\end{proposition}
\begin{proof}
For $i\in V(T_1)$ and $a\in V(T_2)$, let $a_i$ denote the vertex in $U_i$ corresponding to $a$; for $j\in V(T_1)$ and $b\in V(T_2)$, define $b_j$ analogously. 
Observe that $a_i\preceq b_j$ in $T_1\boxtimes T_2$ if and only if one of the two following conditions occurs: $i=j$ and $a\preceq b$ in $T_2$, or $i\neq j$, $i\preceq j$ in $T_1$, and $b=r_2$. By suitably ordering the vertices of $T_1\boxtimes T_2$, this yields the desired formula. 
\end{proof}
\begin{proposition}
\label{expr_bottleneck_matrix_rooted_product}
Let $T_1$ and $T_2$ be rooted trees. Then the bottleneck matrix of $T_1\boxtimes T_2$ is permutationally similar to
\begin{align*}
I_{n_1}\otimes M_2+(M_1-I_{n_1})\otimes J_{n_2}.
\end{align*}
\end{proposition}
\begin{proof}
Using Proposition \ref{expr_path_matrix_rooted_product}, we find that the bottleneck matrix of $T_1\boxtimes T_2$ is permutationally similar to
\begin{align*}
&(I_{n_1}\otimes N_2 +(N_1-I_{n_1})\otimes e_{r_2}e^T)^T (I_{n_1}\otimes N_2 +(N_1-I_{n_1})\otimes e_{r_2}e^T)\\
=&
(I_{n_1}\otimes N_2^T +(N_1^T-I_{n_1})\otimes ee_{r_2}^T) (I_{n_1}\otimes N_2 +(N_1-I_{n_1})\otimes e_{r_2}e^T)\\
=&
I_{n_1}\otimes M_2+
(N_1-I_{n_1})\otimes J_{n_2}
+(N_1^T-I_{n_1})\otimes J_{n_2}\\
&+
(M_1-N_1^T-N_1+I_{n_1})\otimes J_{n_2}\\
=&I_{n_1}\otimes M_2+(M_1-I_{n_1})\otimes J_{n_2}
\end{align*}
as desired.
\end{proof}
\begin{proposition}
\label{expr_neckbottle_matrix_rooted_product}
Let $T_1$ and $T_2$ be rooted trees. Then the neckbottle matrix of $T_1\boxtimes T_2$ is permutationally similar to
\begin{align*}
&I_{n_1}\otimes Q_2+
(N_1^T-I_{n_1})\otimes N_2ee_{r_2}^T+(N_1-I_{n_1})\otimes e_{r_2}e^TN_2^T\\
&+
n_2(Q_1-N_1-N_1^T+I_{n_1})\otimes e_{r_2}e_{r_2}^T.
\end{align*}
\end{proposition}
\begin{proof}
Using Proposition \ref{expr_path_matrix_rooted_product}, we find that the neckbottle matrix of $T_1\boxtimes T_2$ is permutationally similar to
\begin{align*}
&(I_{n_1}\otimes N_2 +(N_1-I_{n_1})\otimes e_{r_2}e^T) (I_{n_1}\otimes N_2 +(N_1-I_{n_1})\otimes e_{r_2}e^T)^T\\
=&
(I_{n_1}\otimes N_2 +(N_1-I_{n_1})\otimes e_{r_2}e^T)
(I_{n_1}\otimes N_2^T +(N_1^T-I_{n_1})\otimes ee_{r_2}^T) \\
=&
I_{n_1}\otimes Q_2+
(N_1^T-I_{n_1})\otimes N_2ee_{r_2}^T+(N_1-I_{n_1})\otimes e_{r_2}e^TN_2^T\\
&+
(Q_1-N_1-N_1^T+I_{n_1})\otimes n_2e_{r_2}e_{r_2}^T
\end{align*}
as desired.
\end{proof}
We now give some lower and upper bounds for the Perron value of the rooted product of two rooted trees $T_1$ and $T_2$. The lower bounds are sharp, meaning that for certain choices of $T_1$ and $T_2$ they are met with equality. The second lower bound involves the Perron entropy introduced in Section \ref{sec_entropy}.
\begin{proposition}
\label{bounds_perron_value_rooted_product}
Let $T_1$ and $T_2$ be rooted trees. Then
\[
\begin{array}{cl}
(i) & \rho(T_1\boxtimes T_2)\geq n_2\rho(T_1);\\[3pt]
(ii) & \rho(T_1\boxtimes T_2)\geq \rho(T_2)+(\rho(T_1)-1)H(T_2);\\[3pt]
(iii) & \rho(T_1\boxtimes T_2)< n_2\rho(T_1)+\rho(T_2).
\end{array}
\]
Moreover, the bounds $(i)$ and $(ii)$ are sharp.
\end{proposition}
\begin{proof}
Using Proposition \ref{expr_bottleneck_matrix_rooted_product}, we have that, for a suitable vertex ordering, the bottleneck matrix of $T_1\boxtimes T_2$ is
\begin{align*}
M=I_{n_1}\otimes M_2+(M_1-I_{n_1})\otimes J_{n_2}=M_1\otimes J_{n_2}+I_{n_1}\otimes(M_2-J_{n_2}),
\end{align*}
so that $M\geq M_1\otimes J_{n_2}$ entrywise.
Hence, by \cite[Corollary~8.1.19]{HJ} and \cite[Theorem~4.2.12]{topics_matrix_analysis},
\begin{align*}
\rho(M)\geq \rho(M_1\otimes J_{n_2})=\rho(M_1)\rho(J_{n_2})=n_2\rho(T_1),
\end{align*}
which proves $(i)$. Also, using the triangle inequality for the spectral norm and \cite[8.2.P15]{HJ}, we obtain
\begin{align*}
\rho(M)&\leq \rho(M_1\otimes J_{n_2})+\rho(I_{n_1}\otimes(M_2-J_{n_2}))\\
&=n_2\rho(T_1)+\rho(M_2-J_{n_2})\\
&< n_2\rho(T_1)+\rho(T_2),
\end{align*}
which proves $(iii)$. Let now $w_1$ and $w_2$ be Perron vectors of norm $1$ for $M_1$ and $M_2$, respectively, and notice that $\|w_1\otimes w_2\|=1$. We find
\begin{align*}
\rho(M)&=\rho(M_1\otimes J_{n_2}+I_{n_1}\otimes(M_2-J_{n_2}))\\
&\geq (w_1\otimes w_2)^T(M_1\otimes J_{n_2}+I_{n_1}\otimes(M_2-J_{n_2}))(w_1\otimes w_2)\\
&=w_1^TM_1w_1w_2^TJ_{n_2}w_2+w_1^TI_{n_1}w_1w_2^T(M_2-J_{n_2})w_2\\
&=\rho(T_1)(e^Tw_2)^2+\rho(T_2)-(e^Tw_2)^2\\
&=\rho(T_2)+(\rho(T_1)-1)H(T_2),
\end{align*}
which proves $(ii)$. 

If $T_2=\mathcal{E}$, then $T_1\boxtimes T_2=T_1$ and $\rho(T_2)=H(T_2)=1$. Hence, the bounds $(i)$ and $(ii)$ hold with equality in this case.
\end{proof}
\noindent Note that the two lower bounds in Proposition \ref{bounds_perron_value_rooted_product} (parts $(i)$ and $(ii)$) are incomparable in general. For instance, letting $T_1=\mathcal{P}_6$ and $T_2=\mathcal{S}_3$ and using \eqref{expr_rho_star_2008}, \eqref{expr_rho_path_2008}, and Proposition \ref{prop_entropy_star_2008}, we find that the first and second bounds are approximately $51.621$ and $51.435$, respectively; letting $T_1=\mathcal{P}_6$ and $T_2=\mathcal{S}_4$, however, we find that the first and second bounds are approximately $68.827$ and $69.035$, respectively.

Given a rooted tree $T$ with root $r$ and order $n$, we let its \textit{root-transmission} $t(T)$ be the sum of the distances of all vertices of $T$ from $r$:
\begin{align*}
t(T)=\sum_{i=1}^n \dist(i,r).
\end{align*}
Observe that the root-transmission of $T$ may be expressed in terms of the sum of the entries of its path matrix $N$:
\begin{align}
\label{expr_root_transmission_path_2008_1735}
t(T)=\sum_{j=1}^n(\dist(j,r)+1)-n=\sum_{i=1}^n\sum_{j\preceq i}1-n=e^TNe-n.
\end{align} 
Also, letting $\gamma=(\gamma_i)$ be the \textit{degree vector} containing the degrees of the vertices in $T$, notice that the moment of $T$ may be expressed as follows:
\begin{align}
\label{expr_moment_path_2008_1727}
\mu(T)&=\sum_{j=1}^n\dist(j,r)\gamma_j=\sum_{j=1}^n(\dist(j,r)+1)\gamma_j-2n+2\notag\\
&=\sum_{i=1}^n\sum_{j\preceq i}\gamma_j-2n+2=\sum_{i=1}^n(N\gamma)_i-2n+2=e^TN\gamma-2n+2.
\end{align}
We are ready to give a result on the moment of the rooted product of two rooted trees.
\begin{proposition}
\label{moment_rooted_product}
Let $T_1$ and $T_2$ be rooted trees. Then
\begin{align*}
\mu(T_1\boxtimes T_2)=\mu(T_1)+n_1\mu(T_2)+2(n_2-1)t(T_1).
\end{align*}
\end{proposition}
\begin{proof}
Let $T=T_1\boxtimes T_2$, and order the vertices of $T$ according to Proposition \ref{expr_path_matrix_rooted_product}. In this proof, given a rooted tree $U$, we shall denote the path matrix, degree vector, and all-ones vector of dimension $|V(U)|$ by $N_U$, $\gamma_U$, and $e_U$, respectively. From \eqref{expr_moment_path_2008_1727}, we have that
\begin{align*}
\mu(T)&=e_T^TN_T\gamma_T-2n_1n_2+2.
\intertext{Also, observe that}
e_T&=e_{T_1}\otimes e_{T_2},\\
\gamma_T&=e_{T_1}\otimes \gamma_{T_2}+\gamma_{T_1}\otimes e_{r_2},\\
N_T&=I_{n_1}\otimes N_{T_2}+(N_{T_1}-I_{n_1})\otimes e_{r_2}e_{T_2}^T,
\end{align*}
where, for the last equation, we have used Proposition \ref{expr_path_matrix_rooted_product}. Using \eqref{expr_root_transmission_path_2008_1735}, we obtain
\begin{align*}
\mu(T)=&(e_{T_1}^T\otimes e_{T_2}^T)(I_{n_1}\otimes N_{T_2}+(N_{T_1}-I_{n_1})\otimes e_{r_2}e_{T_2}^T)(e_{T_1}\otimes \gamma_{T_2}+\gamma_{T_1}\otimes e_{r_2})\\
&-2n_1n_2+2\\
=&
e_{T_1}^TI_{n_1}e_{T_1}e_{T_2}^TN_{T_2}\gamma_{T_2}+e_{T_1}^T(N_{T_1}-I_{n_1})e_{T_1}e_{T_2}^Te_{r_2}e_{T_2}^T\gamma_{T_2}\\
&+
e_{T_1}^TI_{n_1}\gamma_{T_1}e_{T_2}^TN_{T_2}e_{r_2}+e_{T_1}^T(N_{T_1}-I_{n_1})\gamma_{T_1}e_{T_2}^Te_{r_2}e_{T_2}^Te_{r_2}-2n_1n_2+2\\
=&
n_1(\mu(T_2)+2n_2-2)+(t(T_1)+n_1-n_1)(2n_2-2)\\
&+2n_1-2+\mu(T_1)+2n_1-2-(2n_1-2)-2n_1n_2+2\\
=&\mu(T_1)+n_1\mu(T_2)+2(n_2-1)t(T_1)
\end{align*}
as wanted.
\end{proof}

We now consider a third operation on the set of rooted trees. The properties of the corresponding Perron value and moment will be used in Section \ref{sec_Perron_value_and_moment}.  
For a rooted tree $T$ of order $n$, we recursively define its \textit{rooted powers} $T^{\boxtimes k}$ $(k\in\mathbb{N})$ as follows:
\begin{align*}
T^{\boxtimes 0}&\coloneqq \mathcal{E};\\
T^{\boxtimes k}&\coloneqq T\boxtimes T^{\boxtimes k-1}\hspace{1.5cm} (k\geq 1).
\end{align*}
An example is shown in Figure \ref{fig_examples_operations_sum_prod_pow}~$c$. Observe that the order of $T^{\boxtimes k}$ is
\begin{align}
\label{order_T_power_k}
|V(T^{\boxtimes k})|=n^k.
\end{align} 
\begin{proposition}
\label{rho_powers}
Let $T$ be a nontrivial rooted tree and let $k\in\mathbb{N}_{>0}$. Then
\begin{align*}
\rho(T)n^{k-1}\leq \rho(T^{\boxtimes k})\leq \rho(T)\frac{n^k-1}{n-1}. 
\end{align*}
\end{proposition}
\begin{proof}
The first inequality follows directly from part $(i)$ of Proposition \ref{bounds_perron_value_rooted_product}. To prove the second inequality, we use induction on $k$. If $k=1$, the claim is true since $T^{\boxtimes 1}=T\boxtimes \mathcal{E}=T$. If $k\geq 2$, using part $(iii)$ of Proposition \ref{bounds_perron_value_rooted_product} and the inductive hypothesis, we find
\begin{align*}
\rho(T^{\boxtimes k})&=\rho(T\boxtimes T^{\boxtimes k-1})<n^{k-1}\rho(T)+\rho(T^{\boxtimes k-1})\\
&\leq 
n^{k-1}\rho(T)+\rho(T)\frac{n^{k-1}-1}{n-1}=\rho(T)\frac{n^k-1}{n-1},
\end{align*}
which validates the inductive step.
\end{proof}
\begin{proposition}
\label{mu_powers}
Let $T$ be a nontrivial rooted tree and let $k\in \mathbb{N}$. Then
\begin{align*}
\mu(T^{\boxtimes k})=(\mu(T)-2t(T))\frac{n^k-1}{n-1}+2t(T)kn^{k-1}.
\end{align*}
\end{proposition}
\begin{proof}
We use induction on $k$. If $k=0$, the claim is trivially true. If $k\geq 1$, using Proposition \ref{moment_rooted_product} and the inductive hypothesis, we find
\begin{align*}
\mu(T^{\boxtimes k})&=\mu(T\boxtimes T^{\boxtimes k-1})=\mu(T)+n\mu(T^{\boxtimes k-1})+2(n^{k-1}-1)t(T)\\
&=\mu(T)+n\left((\mu(T)-2t(T))\frac{n^{k-1}-1}{n-1}+2t(T)(k-1)n^{k-2}\right)\\
&\quad+2(n^{k-1}-1)t(T)\\
&=(\mu(T)-2t(T))\frac{n^k-1}{n-1}+2t(T)kn^{k-1},
\end{align*} 
thus validating the inductive step.
\end{proof}
\section{Perron value and moment}
\label{sec_Perron_value_and_moment}
The Perron value and the moment can be viewed as two different weights for a rooted tree. In this section we investigate the relation between them. We start off with two examples, concerning rooted stars and paths.
\begin{example}
\label{example_mu_rho_star}
From the expressions for the Perron value and the moment of the rooted star $\mc{S}_n$ reported in \eqref{expr_rho_star_2008} and Example \ref{example_star_as_rooted_sum}, respectively, we observe that
$\rho(\mathcal{S}_n)\sim \mu(\mathcal{S}_n)+2$ as $n$ approaches infinity. In particular,
\begin{align*}
\lim_{n\rightarrow\infty}\frac{\mu(\mathcal{S}_n)}{\rho(\mathcal{S}_n)}=1.
\end{align*}
\end{example}  
\begin{example}
\label{example_mu_rho_path}
Applying the definition \eqref{def_moment}, we find that the moment of the rooted path is $\mu(\mathcal{P}_n)=(n-1)^2$. Using \eqref{expr_rho_path_2008}, we obtain
\begin{align*}
\lim_{n\rightarrow\infty}\frac{\mu(\mathcal{P}_n)}{\rho(\mathcal{P}_n)}&=\lim_{n\rightarrow\infty} 2(n-1)^2\left(1-\cos\left(\frac{\pi}{2n+1}\right)\right)=\frac{\pi^2}{4}\approx 2.47.
\end{align*}
\end{example}
The remaining part of this section is dedicated to inequalities involving the moment and the Perron value of rooted trees. In particular, we will show that $\mu(T)$ is \enquote{almost} an upper bound for $\rho(T)$ (Theorem \ref{mu_rho_f}) and the ratio of these two quantities grows at most linearly in the order of $T$ (Theorem \ref{mu_rho_4_7}) but is not bounded above (Theorem \ref{thm_mu_rho_unbounded}).
\begin{theorem}
\label{mu_rho_f}
Let $f:\mathbb{N}_{>0}\rightarrow\mathbb{R}$ be defined by $f(p)=\frac{1}{2}\left(\sqrt{p^2+2p-3}-p+3\right)$. Then
\begin{align}
\label{eq_mu_rho_f}
\mu(T)\geq \rho(T)-f(n)
\end{align}
for every rooted tree $T$ of order $n$, with equality if and only if $T=\mc{S}_n$.
\end{theorem}
\begin{proof}
We prove the claim by induction on the order of $T$. If $n=1$, then $T=\mathcal{S}_1$, $\mu(T)=0$, $\rho(T)=1$, $f(n)=1$, and, therefore, \eqref{eq_mu_rho_f} holds with equality. Suppose now that the claim holds for rooted trees of order up to $n-1$, and let $T$ be a rooted tree of order $n$. If $T=\mathcal{S}_{n}$, we use the formulae in \eqref{expr_rho_star_2008} and Example \ref{example_star_as_rooted_sum} to check that \eqref{eq_mu_rho_f} holds with equality. If $T\neq \mathcal{S}_n$, let $v$ be a vertex such that $x\coloneqq\dist(v,r)=\max_{w\in V(T)}\dist(w,r)$, where $r$ is the root of $T$. Clearly, $v$ is a pendent vertex and, since $T\neq \mathcal{S}_n$, $x\geq 2$. Consider the rooted tree $\tilde T$ obtained from $T$ by removing $v$ and the unique edge incident with $v$, and having root $r$. The moment of $\tilde T$ is
\begin{align*}
\mu(\tilde{T})=\mu(T)-2x+1.
\end{align*}
Letting $M$ and $\tilde{M}$ denote the bottleneck matrices of $T$ and $\tilde{T}$, respectively, we have that
\begin{align*}
M=P^{-1}\begin{bmatrix}
\tilde M & a \\ 
a^T & x+1
\end{bmatrix} P
\end{align*}
for some $n\times n$ permutation matrix $P$ and some vector $a\in\mathbb{R}^{n-1}$. By using Proposition \ref{prop_fundamental_topics}, we find that
\begin{align*}
\rho(T)\leq \rho(\tilde T)+x+1. 
\end{align*}
Hence, applying the inductive hypothesis to $\tilde T$, we obtain
\begin{align*}
\mu(T)&=\mu(\tilde{T})+2x-1\geq \rho(\tilde T)-f(n-1)+2x-1\\
&\geq \rho(T)-x-1-f(n-1)+2x-1=\rho(T)+x-2-f(n-1)\\
&\geq \rho(T)-f(n-1)>\rho(T)-f(n),
\end{align*}
where the last inequality follows from the fact that $f$ is strictly increasing on $\mathbb{N}_{>0}$.
\end{proof}
\begin{corollary}
For every rooted tree $T$ we have
\begin{align*}
\mu(T)>\rho(T)-2.
\end{align*}
\end{corollary}
Before showing that the ratio of the moment and the Perron value of a rooted tree is at most linear in the number of vertices, we need to prove two technical results.
\begin{lemma}
\label{lemma_N_product_degree_e}
Let $T$ be a rooted tree having root $r$ and let $N$ and $\gamma=(\gamma_i)$ denote its path matrix and degree vector, respectively. Then
\begin{align*}
N(2e-\gamma)=e+e_r.
\end{align*}
\end{lemma} 
\begin{proof}
Using Proposition \ref{inverse_path_matrix}, we see that the claim is equivalent to
\begin{align}
\label{expr_1604_1058}
N^{-1}e+N^{-1}e_r=2e-\gamma.
\end{align}
Observe that $N^{-1}e_r=e_r$. Moreover, $(N^{-1}e)_r=1-\gamma_r$ and $(N^{-1}e)_i=1-(\gamma_i-1)=2-\gamma_i$ if $i\neq r$, so that $N^{-1}e=2e-\gamma-e_r$. From this, \eqref{expr_1604_1058} easily follows.
\end{proof}
Following \cite{Kemeny_bounds}, for $x\in \mathbb{N}_{>0}$ and $y\in \mathbb{N}$, we define the \textit{rooted broom} $\mc B(x,y)$ as the rooted tree obtained by attaching $y$ pendent vertices to an endpoint of a path of $x$ vertices and by letting the other endpoint be the root (if $x=1$, we let $\mc B(x,y)\coloneqq\mc{S}_{y+1}$). We also let $\mc B(0,1)\coloneqq\mc{E}$.
\begin{proposition}
\label{max_mu_fixed_order}
Let $T$ be a rooted tree of order $n$. Then
\begin{align*}
\mu(T)\leq (n-1)^2
\end{align*}
with equality if and only if $T=\mathcal{P}_n$.
\end{proposition}
\begin{proof}
Let $\delta$ denote the diameter of $T$. Using \cite[Proposition~5.2]{Kemeny_bounds}, we have that
\begin{align*}
\mu(T)\leq 2n\delta-\delta^2-n-\delta+1\eqqcolon p(n,\delta),
\end{align*}
with equality if and only if $T=\mc B(\delta,n-\delta)$. Observing that, for $\delta\leq n-1$,
\begin{align*}
\frac{\partial p}{\partial \delta}=2n-2\delta-1\geq 2\delta+2-2\delta-1=1>0,
\end{align*}
we deduce that 
\begin{align*}
p(n,\delta)\leq p(n,n-1)=(n-1)^2,
\end{align*}
with equality if and only if $\delta= n-1$. We conclude that $\mu(T)\leq (n-1)^2$, with equality if and only if $T=\mc B(n-1,1)=\mathcal{P}_n$. 
\end{proof}
\begin{theorem}
\label{mu_rho_4_7}
Let $T$ be a rooted tree of order $n$. Then
\begin{align*}
\mu(T)< \frac{4}{7}n\rho(T).
\end{align*}
\end{theorem}
\begin{proof}
If $n=1$ the claim is trivial, so assume that $n\geq 2$. Let $r$, $N$, $M$, and $\gamma=(\gamma_i)$ denote the root, path matrix, bottleneck matrix, and degree vector of $T$, respectively. 
Using \eqref{expr_moment_path_2008_1727}, Lemma \ref{lemma_N_product_degree_e}, and the identity $N^Te_r=e$, we find
\begin{align*}
4e^TMe&=(2Ne)^T(2Ne)=(\gamma^TN^T+e^T+e_r^T)(N\gamma+e+e_r)\\
&=\gamma^TN^TN\gamma+e^Te+e_r^Te_r+2e^TN\gamma+2e_r^TN\gamma+2e^Te_r\\
&=\gamma^TM\gamma+n+1+2(\mu(T)+2n-2)+2(2n-2)+2\\
&=2\mu(T)+\gamma^TM\gamma+9n-5.
\end{align*} 
Observe that, entrywise,
$M\geq ee^T+\Delta$,
where $\Delta$ is the $n\times n$ diagonal matrix whose $(i,i)$'th entry is $\dist(i,r)$. We obtain
\begin{align*}
\gamma^TM\gamma&\geq \gamma^Tee^T\gamma+\gamma^T\Delta\gamma
\geq \gamma^Tee^T\gamma+e^T\Delta\gamma
=(2n-2)^2+\sum_{i=1}^n\dist(i,r)\gamma_i\\
&=4(n-1)^2+\mu(T)
\geq 5\mu(T),
\end{align*}
where, for the last inequality, we have applied Proposition \ref{max_mu_fixed_order}. We conclude that
\begin{align*}
4e^TMe &\geq 7\mu(T)+9n-5> 7\mu(T)
\end{align*}
and, hence,
\begin{align*}
\mu(T)<\frac{4}{7}e^TMe=\frac{4}{7}n\frac{e^TMe}{e^Te}\leq \frac{4}{7}n\rho(T)
\end{align*}
as desired.
\end{proof}
Example \ref{example_mu_rho_star} and Example \ref{example_mu_rho_path} show that for both the rooted star and the rooted path -- which may be considered the two \textit{extremal} classes of rooted trees in many respects -- the ratio of the moment and the Perron value is bounded above. This could suggest trying to sharpen Theorem \ref{mu_rho_4_7}, to show that the ratio of the moment and the Perron value of generic rooted trees is bounded above. However, it turns out that certain families of rooted trees exhibiting a \enquote{fractal} structure provide a counterexample to this intuition. We now use the results in Section \ref{sec_rooted_sums_product_powers} to identify two such families for which the $\mu/\rho$ ratio is unbounded. As a consequence, we shall prove the next result.
\begin{theorem}
\label{thm_mu_rho_unbounded}
For any $\alpha\in\mathbb{R}$ there exists a rooted tree $T$ such that $\mu(T)>\alpha\,\rho(T)$.
\end{theorem}
Consider first the class of rooted Bethe trees defined in \eqref{def_bethe_trees}. Combining the results in Proposition \ref{rho_ennary_tree} and
Proposition \ref{moment_ennary_tree}, we see that, for fixed $k\geq 2$ and large $p$,
\begin{align*}
\frac{\mu(\mc B_{p,k})}{\rho(\mc B_{p,k})}\geq\frac{2pk^{p+1}-3k^{p+1}-2pk^{p}+k^{p}+k^2+k}{k^{p+1}-pk-k+p}\sim 
\frac{1}{k}-3
+2\left(1-\frac{1}{k}\right)p
\end{align*}
and, in particular,
\begin{align}
\label{unbounded_ratio_mu_rho_bethe}
\lim_{p\rightarrow\infty}\frac{\mu(\mc B_{p,k})}{\rho(\mc B_{p,k})}=\infty.
\end{align}  
\begin{observation}
\label{remark_log_bethe}
From \eqref{order_bethe_trees}, notice that $p=\mc{O}(\ln(|V(\mc{B}_{p,k})|))$ as $p\rightarrow\infty$.
\end{observation}
Let now $T$ be a nontrivial rooted tree of order $n$. Using Proposition \ref{rho_powers} and Proposition \ref{mu_powers}, we see that, for large $k$, 
\begin{align*}
\frac{\mu(T^{\boxtimes k})}{\rho(T^{\boxtimes k})}
&\geq 
\left((\mu(T)-2t(T))\frac{n^k-1}{n-1}+2t(T)kn^{k-1}\right)\frac{n-1}{\rho(T)(n^k-1)}\\
&=\frac{\mu(T)-2t(T)}{\rho(T)}+2\frac{t(T)}{\rho(T)}k\frac{n^k-n^{k-1}}{n^k-1}\\
&\sim
\frac{\mu(T)-2t(T)}{\rho(T)}
+2\frac{t(T)}{\rho(T)}\left(1-\frac{1}{n}\right)k
\end{align*}
and, in particular, 
\begin{align}
\label{unbounded_ratio_mu_rho_powers}
\lim_{k\rightarrow \infty}\frac{\mu(T^{\boxtimes k})}{\rho(T^{\boxtimes k})}=\infty.
\end{align}
\begin{observation}
\label{remark_log_powers}
From \eqref{order_T_power_k}, notice that $k=\mc{O}(\ln(|V(T^{\boxtimes k})|))$ as $k\rightarrow\infty$.
\end{observation}
Using either \eqref{unbounded_ratio_mu_rho_bethe} or \eqref{unbounded_ratio_mu_rho_powers}, one proves Theorem \ref{thm_mu_rho_unbounded}. In light of Observation \ref{remark_log_bethe} and Observation \ref{remark_log_powers}, we conclude with the following conjecture.
\begin{conjecture}
There exists $\alpha_0\in\mathbb{R}$ such that $\mu(T)\leq \alpha_0 \ln(n)\rho(T)$ for any rooted tree $T$ of order $n$.
\end{conjecture}

\bigskip
\section*{Acknowledgements}
Research supported by a Doctoral Research Fellowship at the Faculty of Mathematics and Natural Sciences, University of Oslo.
The author is grateful to Enide Andrade and Geir Dahl for many fruitful discussions and important comments on this work.
\bigskip
%


\bibliographystyle{plain}





\end{document}